\documentclass[12pt,oneside]{article}
\usepackage{amssymb,amsmath,amstext,amsthm,amsfonts}
\usepackage[ansinew]{inputenc} 
\usepackage{graphicx}
\usepackage[mathscr]{eucal}
\usepackage[pdfnewwindow=true]{hyperref}

\newcommand{\R}{\mathbb{R}}
\newcommand{\T}{\mathbb{T}}
\newcommand{\Q}{\mathbb{Q}}
\newcommand{\N}{\mathbb{N}}
\newcommand{\Z}{\mathbb{Z}}

\newcommand{\NN}{\mathbb{N}}
\newcommand{\Escr}{\mathscr{E}}
\newcommand{\Fscr}{\mathscr{F}}

\newcommand{\Mod}[1]{\left\vert{#1}\right\vert}

\newcommand{\nrm}[1]{\left\|#1\right\|}

\newcommand{\spnrm}[1]{\vert\!\vert\!\vert{#1}\vert\!\vert\!\vert}

\newcommand{\Lin}[3]{ \mathscr{L}^{#1}({#2},{#3}) }
\newcommand{\Linsim}[3]{ \mathscr{L}^{#1}_{sym}({#2},{#3}) }
\newcommand{\Hom}[3]{ \mathscr{H}^{#1}({#2},{#3}) }
\newcommand{\Pol}[3]{ \mathscr{P}^{#1}({#2},{#3}) }

\newcommand{\diag}{\mbox{diag}}
\newcommand{\coment}[1]{}

\newcommand{\CD}[3]{ \mathscr{C}^{#1}( {#2},{#3}) }

\newcommand{\ovl}{\overline}

\newcommand{\slim}[1]{\stackrel{\longrightarrow}\lim_{{#1}\rightarrow\infty} }
\newcommand{\bin}[2]{\left(^{#1}_{#2}\right)}
\newtheorem{lema}{Lema}
\newtheorem{prop}{Proposition}
\newtheorem{teor}{Theorem}

\newcommand{\dem}{ \par\medbreak\noindent{\bf
Proof. }\enspace} 
\newcommand{\demof}[1]{ \par\medbreak\noindent{\bf
{#1}. }\enspace} 
\newcommand{\cqd}{\hfill
$\sqcup\!\!\!\!\sqcap\bigskip$}

\title{ Extension of functions with bounded finite differences}
\author{ P. Duarte$\,^{1}$,  \; M. J. Torres$\,^{2}$ }
\date{}

\begin{document}
\maketitle

{\small
\noindent
$^1$ {\it CMAF, Departamento de Matem\'{a}tica, Faculdade de Ci\^{e}ncias da Universidade de Lisboa, Portugal
(e-mail: pedromiguel.duarte@gmail.com).
Supported by Funda\c{c}\~{a}o para a Ci\^{e}ncia e a Tecnologia, Financiamento Base 2008 - ISFL/1/209.} \\
\noindent
$^2$ {\it CMAT, Departamento de Matemática, Universidade do Minho, Campus de Gualtar, 4700-057 Braga,
Portugal (e-mail: jtorres@math.uminho.pt). }
}

\bigskip

\noindent
{\bf Abstract }
{\em  We prove that functions  defined on a lattice in the torus $\T^d$
with bounded finite differences can be smoothly extended to the whole torus,
and relate the bounds on the extension's derivatives with bounds on the
original function's finite differences.
}

 \bigskip

\section{Introduction}

In 1934, Whitney posed the problem of how to recognize whether a function $f$
defined on a closed subset $X$ of $\R^n$ is the restriction of a function
of class $C^m$. Whitney himself solved the one-dimensional case (i.e., for $n=1$)
in terms of finite differences~\cite{W1, W2, W3}, giving the classical Whitney's extension theorem. 
A geometrical solution for the case $C^1(\R^n)$  was given by G. Glaeser ~\cite{G},
who introduced 
 a geometric object called the ``iterated paratangent space''. Glaeser's 
paper influenced all the later work on Whitney's problem. A variant of 
 Whitney's problem replaces $C^m(\R^n)$ by $C^{m,\omega}(\R^n)$, the 
space of $C^m$ functions whose $m^{\rm th}$ derivative have a given modulus of
continuity $\omega$. The problem for $C^{m,\omega}(\R^n)$ is well-understood
due to the work of Brudnyi and Shvartsman~\cite{B}, \cite{BS1,BS2,BS3,BS4}, 
\cite{S1,S2,S3} and Fefferman~\cite{F1}.
The correct notion of an iterated paratangent bundle, relevant for 
$C^m(\R^n)$, was introduced by Bierstone-Milman-Pawlucki~\cite{BMP1}
who proved an extension theorem for subanalytic sets. In~\cite{BMP1}, 
Bierstone-Milman-Pawlucki introduced a necessary geometric criterion,
for $C^m(\R^n)$, involving limits of finite differences, and conjectured 
that this criterion is sufficient at least if $X$ has a ``tame topology''.
In~\cite{F2} C. Fefferman solved the $C^m(\R^n)$ 
Whitney's problem. After that Bierstone-Milman-Pawlucki~\cite{BMP2} verified that
the conjectures of~\cite{BMP1} with the paratangent bundle there
replaced by a natural variant are equivalent to Fefferman's solution of Whitney's problem~\cite{F2}.

We base our work here in the extension on the classical
Whitney's extension theorem proved by C. Fefferman in~\cite{F1}.
Consider any set  $\Gamma\subset \R^d$, possibly finite.
We shall say that a function $f:\Gamma\to \R$ satisfies
the {\em $C^k$-Whitney extension condition with constant $M\geq 0$} \, iff\,
there is a family of polynomials $\{P_x\}_{x\in\Gamma}$ such that
for every $x,y\in\Gamma$, every $m=0,1,\ldots, k-1$, and 
every multi-index $(i_1,\ldots, i_m)\in\N^m$,
\begin{enumerate}
\item $P_x(x)=f(x)$  
\item $\Mod{ D^m_{( e_{i_1},\ldots, e_{i_m}) } (P_x)(x)}\leq M$
\item $\Mod{ D^m_{( e_{i_1},\ldots, e_{i_m}) } (P_x-P_y)(x)}\leq M\,\nrm{x-y}^{k-m}$
\end{enumerate}
where $D^m_{( e_{i_1},\ldots, e_{i_m})}$ stands for $\partial^m/ \partial x_{i_1}\cdots \partial x_{i_m}$.
Fefferman's theorem can be stated as follows:

\bigskip

\noindent
{\bf Theorem }(C. Fefferman) {\em
Given $d,k\in\N$ there are constants $A=A(d,k)>0$ and $n=n(k,d)\in\N$ such that
for any $M\geq 0$ and any function $f:\Gamma\to \R$, defined on a finite set $\Gamma\subset \R^d$, 
if the restriction $f\vert_E:E\to \R$ to any subset of $E\subseteq \Gamma$ with $\# E\leq n$ 
satisfies the $C^k$-Whitney extension condition with constant $M$, then
$f$ has an extension of class $C^k$ $F:\R^d\to \R$ with 
$\nrm{F}_{C^{k}(\R^d)} \leq A\,M$.
}

\bigskip

See Theorem B of~\cite{F1}.
Although this is a strong result, the Whitney's condition is not
easy to check in practice. We prove here, in the very particular setting
where $\Gamma$ is a lattice in the torus $\T^d$,
that Whitney's condition on $f$ is implied by the assumption
that $f$ has bounded finite differences of order $k$.
Given $m\in\Z_+^d$ consider the lattice
$$\Gamma_m= ( {m_1}^{-1}\Z \times \ldots \times  {m_d}^{-1}\Z ) /\Z^d\subset \T^d \;. $$
We denote by $\spnrm{\Delta^k f}_{\Gamma_m}$ the least upper bound of all
differences of order $k$ of function $f$ over ${\Gamma_m}$, and by
$\spnrm{D^k F}_{\T^d}$ the least upper bound of all
derivatives of order $k$ of a class $C^k$ function $F:\T^d\to\R$.
With this notation we prove that

\bigskip

\noindent
{\bf Theorem A} {\em
Given $d,k,\ell\in\N$ there is a constant $A=A(d,k,\ell)>0$ such that
for any function $f:{\Gamma_m}\to \R$\, defined on the lattice ${\Gamma_m}\subset \T^d$, \,
where $m\in\Z^d_+$ is such that $m_i/m_j\leq \ell$ for every $i,j=1,\ldots, d$,\, then 
$f$ has an extension of class $C^{k-1}$ $F:\T^d\to \R$ with 	\, 
${\rm Lip}(D^{k-1} F) \leq A\,\spnrm{\Delta^k f}_{\Gamma_m}$.
}

\bigskip

We believe this result should be true with $A=1$, but 
we do not prove it. Of course a direct approach
(not using Fefferman's theorem) is necessary for such purpose.
With the same type of argument we can prove the following theorem, 
which is likely to be known, but whose proof
we are not aware in the literature.

\bigskip

\noindent
{\bf Theorem B} {\em Let $\Omega\subset\R^d$ be an open convex set
and  $f:\Omega\to \R$ a continuous function.
Then $f$ has bounded finite differences of order $k$ \, iff\,
$f$ is of class $C^{k-1}$ with Lipschitz derivative of order $k-1$.
Furthermore,  \, 
$\spnrm{\Delta^k f}_\Omega={\rm Lip}(D^{k-1} f)$.
}

\bigskip

See theorem~\ref{teor:B}.
Theorem A should also hold for some class of domains $\Sigma$ which are
proper subsets of lattices $\Gamma_m$. We would like to characterize
this class of subsets. More precisely, we pose the following problem

\bigskip

\noindent
{\bf Problem} (Extensibility of discrete $L^k$ functions) {\em  
Characterize the pairs $(k,\Sigma)$,
with $\Sigma\subset \Gamma_m$ and $k\in\N$, for which given any function $f:\Sigma\to\R$ 
(resp. $f:\Sigma/\Z^d\to\R$) with
$\spnrm{\Delta^k f}_\Sigma<+\infty$ there is an extension $\tilde{f}:\Gamma_m\to\R$ 
(resp. $\tilde{f}:\Gamma_m/\Z^d\to\R$)
such that \, $\tilde{f}=f$ over $\Sigma$  (resp. $\Sigma/\Z^d$) \, and \, 
$\spnrm{\Delta^k \tilde{f}}_{\Gamma_m}=\spnrm{\Delta^k f}_\Sigma$.
}

\bigskip

\section{Finite difference operators}

Given two normed spaces $X$ and $Y$, let $Y^X$ denote the space of
all functions $f: X\to Y$.
For each vector $u\in X$, we define a {\em difference operator}\;
$\Delta_u :Y^X\rightarrow Y^X$,
$$ \Delta_u f(x) = f(x+u)-f(x) = \left( (\mbox{Id}-\tau_u)\,f\right)(x)\;,$$
where $\tau_u$ represents the shift operator\; $(\tau_u f)(x)= f(x+u)$.
Remark that 
\begin{equation}
\label{abel}
\Delta_u\circ \Delta_v = \Delta_v\circ \Delta_u\;,
\end{equation}
for all vectors $u, v\in X$, since
$$ \tau_v\circ\tau_u
 = \tau_{u+v } = \tau_u\circ\tau_v\;.$$

Another key property of the difference operator $\Delta$ is the following
kind of additivity:
\begin{equation}
\label{aditiv}
 \Delta_{u+v} f(x) = \Delta_u f(x) + \Delta_v f(x+u)\;.
\end{equation}

More generally, given a multi-vector $u=(u_1,\cdots,u_k)\in X^k$, 
we define the  {\em finite difference operator of order $k$}, along the multi-vector $u$, 
as the composition operator :\;
$\Delta^k_u :Y^X \rightarrow Y^X$,\;
$\Delta^k_u = \Delta_{u_1}\circ\Delta_{u_2}\circ\cdots\circ \Delta_{u_k}$.

\bigskip

In order to characterize this operator consider the 
discrete cube  $I^k=\{0,1\}^k$ as 
a set of multi-indices. Denote by $\alpha$
a generic element of the cube $I^k$,
and write $\alpha=(\alpha_1,\cdots,\alpha_k)$, where each
$\alpha_i$ represents a binary digit, $\alpha_i=0$ or $\alpha_i=1$.
The set $I^k$ is partially ordered by the relation
$$\alpha\leq \beta \quad \Leftrightarrow \quad \alpha_i\leq \beta_i,\;
\text{ for all }\; i=1,\ldots,k\;.$$
We shall write
$\Mod{\alpha}=\alpha_1+\ldots +\alpha_k$\; and \;
$\alpha\cdot u = \alpha_1\,u_1+\ldots + \alpha_k\,u_k$.
Using this notation we have

\begin{prop}
\label{Diff}
$$\Delta^k_u = \sum_{\alpha\in I^k} 
(-1)^{k-\Mod{\alpha}}\,\tau_{\alpha\cdot u}\;, $$
In particular

$$\Delta^k_u f (x) = \sum_{\alpha\in I^k} 
(-1)^{k-\Mod{\alpha}}\,f( x+\alpha\cdot u)\;. $$
\end{prop}

\bigskip

The operator $\Delta^k_u$ also  acts
on partial functions. Assume $f\in Y^\Gamma$ is some function
with domain $\Gamma\subseteq X$. Given the multi-vector   
$u\in X^{k}$,\,
we define
$$\Gamma^{u} = \{ \,x\in \Gamma\,:\;
x+\alpha\cdot u\in \Gamma\;\text{ for all }\; 
\alpha\in I^k\;\}\;.$$ 
Then $\Delta^k_u f$ is a new function defined on $\Gamma^{u}$.

\bigskip

The difference operator $\Delta_u$ is the discrete analogous of
the directional derivative operator
$$ D_u f (x) = \lim_{h\to 0} \frac{1}{h}\,\Delta_{h\,u} f(x) \;.$$
Under general conditions, the directional derivative
operators satisfy for all $u, v\in X$,
$$ D_u\circ D_v = D_v\circ D_u\quad \text{and}\quad
 D_{u+v}  = D_u  + D_v \;, $$
which are the infinitesimal 
equivalents to~(\ref{abel}) and~(\ref{aditiv}).
The finite differences $\Delta^k_u f(x)$  are also discrete analogs
of the higher order directional derivatives
$$ D^k_u f (x) = \left(D_{u_1} \circ \ldots \circ D_{u_k} \right) f (x)\;. $$
The operators $D_u$ are well defined
on the space $\CD{\infty}{\Gamma}{Y}$ of smooth functions over
an open set $\Gamma\subset X$.
If $f\in\CD{k}{\Gamma}{Y}$ then
$\Delta^k_u f(x)$ is, in the following sense,  an average of $D^k_u f$ values.

\begin{prop}[Mean Value Theorem]
\label{mean value theorem}
Given a function $f\in\CD{k}{\Gamma}{Y}$,
defined over an open set $\Gamma\subseteq X$,
containing the parallelogram \, $\{ \, x+t\cdot u\in X :\, t\in [0,1]^k\,\}$,
then
$$ \Delta^k_u f(x) = \int_{[0,1]^k}  D^k_u f(x+t\cdot u)\,dt\;.$$
\end{prop}

\dem
For $k=1$, this is the usual mean value theorem.
Assume this proposition holds for $k-1$.
Given $f\in\CD{k}{\Gamma}{Y}$, we have by induction hypothesis,
$$\Delta^{k-1}_u f(x) = \int_{[0,1]^{k-1}}  D^{k-1}_u f(x+t\cdot u)\,dt\;.$$
Therefore
\begin{eqnarray*}
\Delta_v \Delta^{k-1}_u f(x) &=& \Delta_v \int_{[0,1]^{k-1}}  D^{k-1}_u f(x+t\cdot u)\,dt\\
&=&  \int_{[0,1]^{k-1}} \Delta_v (D^{k-1}_u f)(x+t\cdot u)\,dt\\
&=&  \int_{[0,1]^{k-1}} \int_0^1  D_v D^{k-1}_u f(x + t\cdot u + s\,v)\,dt\;,
\end{eqnarray*}
which proves that the same formula holds for $k$.
\cqd

\bigskip

Next we derive an important formula which in some sense is a discrete version
of proposition~\ref{mean value theorem}.
Given a multi-index 
$n=(n_1,n_2, \ldots,n_k)\in\N^k$,\;
we define\; 
$$ [n] =\{ \, m\in\N^k\; :\; 0\leq m_i<n_i\quad\forall\; i=1,\ldots,k\; \}\;,$$
and $\ovl{n}=n_1\,n_2\,\ldots \,n_k$.
Notice that $[n]$ has  $\ovl{n}$ elements.
Given a multi-vector $u\in X^k$,\,
we set $n\cdot u = n_1 u_1+\ldots+n_k u_k\in X$  
and write
$\Mod{n}=n_1+n_2+\ldots +n_k$.

\begin{prop}
Given $f\in Y^\Gamma$ with $\Gamma\subseteq X$, $x\in\Gamma$, $u\in X^k$  
and $n\in\N^k$ such that \,
$x+ j\cdot u \in\Gamma$,\, for all $j\in [n]$, then 
\begin{equation}
\label{genaditiv}
\Delta^k_{n\,u} f(x) = \sum_{j\in [n]} \Delta^k_{u} f\left( x+ j\cdot u \right)\;.
\end{equation}
\end{prop}

\dem
Relation~(\ref{aditiv}) is equivalent to $\Delta_{u+v}=\Delta_u+\tau_u\circ\Delta_v$.
This formula can be generalized to an arbitrary number of terms. 
Setting all them equal to $u_i$ we get 
$\Delta_{n_i\,u_i} = \sum_{j=0}^{n_i-1} 
\tau_{j\cdot u_i} \circ \Delta_{u_i}$.
Relation~(\ref{genaditiv}) follows by combining
these decompositions.
\cqd

\bigskip

\section{Polynomials}
\label{polynomials}

As in the previous section, let $X$ and $Y$ denote normed vector spaces.
Let $\Lin{k}{X}{Y}$ denote the space of continuous $k$-multi-linear maps
$f:X^k= X\!\times\!X\!\times\ldots\!\times\!X\rightarrow Y$. The $k$-multi-linearity 
assumption means linearity in each of the function $k$ arguments.
We shall denote by $\Linsim{k}{\R^d}{Y}$, the subspace of 
continuous $k$-multi-linear symmetric maps. We consider on this space the usual norm
$$ \nrm{f}=\sup \left\{ \, \frac{ \nrm{f(u_1,\ldots, u_k)} }{\nrm{u_1}\cdots \nrm{u_k}} 
\; :\, u_1,\ldots, u_k\in X-\{0\}\, \right\}\;.$$

We define the diagonal operator $\diag : \Linsim{k}{\R^d}{Y}\rightarrow Y^X\;$ by
$$\diag(f) (x) = f(\overbrace{x,x,\ldots,x}^{k\;\text{times}}) \;.$$
Any map in the image of this operator is called a  {\em homogeneous polynomial of degree $k$}.
The space of homogeneous polynomials of degree $k$ is, therefore,
$$\Hom{k}{X}{Y}=\{ \, \diag(f)\; :\; f\in\Linsim{k}{\R^d}{Y}\;\}.$$
Every finite sum of homogeneous polynomials is called a {\em polynomial}.
The degree of a polynomial is the highest degree of all its
homogeneous components. The space of  
polynomials of degree $\leq k$, denoted by $\Pol{k}{X}{Y}$, is the direct sum
$$\Pol{k}{X}{Y}=\bigoplus_{n=0}^k \Hom{n}{X}{Y}\;.$$

%
%

\bigskip

\begin{prop}
\label{polchar}
Given $f\in \CD{0}{X}{Y}$, the following  statements are equivalent:
\begin{enumerate}
\item[(a)] $f$ is a polynomial of degree $\leq k$, i.e., $f\in\Pol{k}{X}{Y}$,
\item[(b)] $\Delta^{k+1}_u f(x)=0$, for all $x\in X$ and  $u\in X^{k+1}$.
\end{enumerate}
\end{prop}

\bigskip
\bigskip

The proof is based on proposition~\ref{derivn} bellow.

\begin{lema}
\label{derivum}
Given $\xi\in\Linsim{k}{\R^d}{Y}$, the homogeneous polynomial
$f=\diag(\xi)$ satisfies
$$ \Delta_u f (x)= k\,\xi(u,x,\ldots,x) + r(x) \;\;
\text{ with }\, r\in\Pol{k-2}{X}{Y}\;. $$
\end{lema}

\dem
Using the binomial Newton formula,
\begin{eqnarray*}
f(x+u) &=&\xi(x+u,x+u,\ldots,x+u )\\
&=& \sum_{i=0}^k  \bin{k}{i}\,\xi( u^{(i)},x^{(k-i)})\\
&=& f(x) +k\,\xi(u,x,\ldots,x)+\cdots
\end{eqnarray*}
where the points refer to a sum of homogeneous polynomials of degree $\leq k-2$.
The notation $u^{(i)}$ refers to a list of $i$
vectors equal to $u$.
\cqd

From this lemma, we get by induction:

\begin{prop}
\label{derivn}
Given $\xi\in\Linsim{k}{\R^d}{Y}$, the homogeneous polynomial
$f=\diag(\xi)$ satisfies
\begin{enumerate}
\item[(1)]
$\displaystyle{ \Delta^k_u f (x)= k!\,\xi(u) }$\;
for every \, $u\in X^k$.
\item[(2)]
$\displaystyle{ \Delta^m_u f (x)= 0 }$\;
for every \, $m\geq k+1$ and \; $u\in X^m$.
\end{enumerate}
\end{prop}

\bigskip
\bigskip

\demof{Proof of Proposition~\ref{polchar}}
That $(a)\Rightarrow (b)$ follows by item (2) of proposition~\ref{derivn}.
The converse is proved by induction in $k$.
By definition, a degree zero polynomial is just 
a constant function. For $k=0$, the condition $(b)$ says that $\Delta_u f(x)=0$
for all $x, u\in X$, which is equivalent to  $f$ being constant.
Assume that condition $(b)$ holds for $k$, and $(b)\Rightarrow (a)$ holds for $k-1$.
Define $\xi:X^k\to Y$ to be
$\xi(u)= \frac{1}{k!}\,\Delta^k_u f(0)$.
Notice that, because of (b),
we also have $\xi(u)= \frac{1}{k!}\,\Delta^k_u f(x)$, for all $x\in X$.
By property~(\ref{abel}) of difference operators,
we get the symmetry of function $\xi$.
Using property~(\ref{aditiv}), given $u_1, u_2\in X$ and $v\in X^{k-1}$,
\begin{eqnarray*}
\xi(u_1+u_2,v) &=& \frac{1}{k!}\,\Delta^{k-1}_v \Delta_{u_1+u_2} f(0)\\
&=& \frac{1}{k!}\,\Delta^{k-1}_v \Delta_{u_1} f(0)+
 \frac{1}{k!}\,\Delta^{k-1}_v \Delta_{u_2} f(u_1)\\
 &=& \xi(u_1,v) + \xi(u_2,v)\;.
\end{eqnarray*}
This proves that $\xi$ is $k$-additive.
By continuity of $f$, $\xi$ is also continuous, and this implies its
$k$-multi-linearity. Therefore,  $\xi\in\Linsim{k}{\R^d}{Y}$.
Take now $g=f-\diag(\xi)$.
By item (1) of proposition~\ref{derivn}
we have $\Delta^k_u g(x)=\Delta^k_u f(x) -k!\,\xi(u)= \Delta^k_u f(x) -\Delta^k_u f(0)=0$,
for all $x\in X$ and $u\in X^k$. Thus, by induction hypothesis,
$g$ is a polynomial of degree $\leq k-1$, and
 then $f= g +\diag(\xi)\in\Pol{k}{X}{Y}$.
 This shows that $(b)\Rightarrow (a)$ holds for $k$,
which ends the induction proof.
\cqd

\bigskip
\bigskip

By item(1) of proposition~\ref{derivn} we have
\begin{lema}
The linear map \; $\diag:\Linsim{k}{\R^d}{Y}\to \Hom{k}{X}{Y}$\;
is an isomorphism, the inverse mapping being such that\;
$\diag^{-1}(f)(u)= \frac{1}{k!}\,\Delta^k_u f(0)$.
\end{lema}

\bigskip

\section{$L^k$-spaces}
\label{Lk:spaces}

Let $G\subseteq X$ be an additive subgroup.
We say that $\Gamma\subseteq X$ is {\em $G$-convex}\, iff\,
for every $x\in X$ and $n\in\N$ such that
$x, x+ n\,u\in\Gamma$ we have $x+i\,u\in\Gamma$ for every $1\leq i \leq n$.
From now on every domain $\Gamma$ will be a $G$-convex w.r.t.
some group additive subgroup $\Gamma\subseteq X$.
In fact we shall only be interested in two cases: either
$G=X$ and $\Gamma$ is convex in the usual sense, or else
$G$ is a lattice in $X$.
We write $G^\ast := G-\{0\}$.

We define the {\em multi-norm} of a multi-vector $u\in X^k$
to be the product
$$\nrm{u}^k = \nrm{u_1}\cdot\nrm{u_2}\cdot \ldots \cdot \nrm{u_k} \;.$$

We say that a  function $f\in Y^\Gamma$ is  {\em Lipschitz of order $k$}
on $\Gamma$  \; if and only if \; $f$ is continuous, and 
there is some constant $C>0$ such that for all $u\in G^k$ \; and \,
$x\in \Gamma^{u}$,\;
$\nrm{\Delta^{k}_u f(x)}\leq C\,\nrm{u}^k$.
We denote by $L^k(\Gamma,Y)$ the space of functions
$f\in Y^\Gamma$ which are Lipschitz of order $k$.
Let us now introduce a norm on this space.
Given $u\in G^k$, we denote by  $\spnrm{\Delta^k_u f}_\Gamma$ the least upper bound 
$$ \spnrm{\Delta^k_u f}_\Gamma = 
\sup_{x\in \Gamma^{u}} \nrm{ \Delta^k_u f(x)}
\;\in [0,+\infty] \;.$$
Remark that $\spnrm{\Delta^k_u f}_\Gamma =0$ when $\Gamma^u=\emptyset$.
We also denote by  $\spnrm{\Delta^k f}_\Gamma$  the least upper bound
$$ \spnrm{\Delta^k f}_\Gamma = 
\sup_{ u\in (G^\ast)^k} \frac{ \spnrm{\Delta^k_{u} f}_\Gamma }{\nrm{u}^k}
\;\in [0,+\infty] \;.$$
The function $f\mapsto \spnrm{\Delta^k f}_\Gamma$ is a pseudo-seminorm on 
$Y^\Gamma$ such that
\[ L^k(\Gamma,Y)=\{\, f\in Y^\Gamma\,:\, \spnrm{\Delta^k f}_\Gamma<+\infty \,\} \;.\]
For $k=0$, $\Delta^0$ is the identity operator and 
$\spnrm{\Delta^0 f}_\Gamma = \sup_{x\in\Gamma} \nrm{ f(x) }$ is the usual
$C^0$-norm.
We shall consider the following norm on the space $L^k(\Gamma,Y)$
\[ \spnrm{f}_{\Gamma, k} =\max\{ \spnrm{\Delta^0 f}_\Gamma , \spnrm{\Delta^k f}_\Gamma \}\;.\]
It can easily be checked that
\begin{prop}
The normed space $\left( L^k(\Gamma,Y),\, \spnrm{\cdot}_{\Gamma, k}\right)$ is a Banach space.
\end{prop}

\bigskip

\begin{prop}[Monotonicity]
\label{monot:prop1} If $\Gamma\subset X$ is $G$-convex, given $f\in Y^\Gamma$  and  $u\in G^{k}$,
\begin{equation}\label{monot}
\spnrm{\Delta^k_{n\,u} f}_\Gamma \leq \ovl{n}\,\spnrm{ \Delta^k_{u}f}_\Gamma\;.
\end{equation}
\end{prop}

\dem
Follows from~(\ref{genaditiv}), where the $G$-convexity of $\Gamma$ is used to ensure that\\
$\displaystyle{ x\in \Gamma^{n u}\quad \Rightarrow \quad x+j\cdot u\in \Gamma,\;
\quad \forall\; j\in [n] }$\;.
\cqd

%
%
%
%

\bigskip

\begin{prop}
For any function $f \in L^{k}(\Gamma,Y)$ and every
 $u, u' \in G^{k}$,
\[ \nrm{ \Delta^k_{u'} f(x) - \Delta^k_{u} f(x) } \leq
\spnrm{\Delta^k f}_\Gamma\,\sum_{i=1}^k
\Mod{u_1'}\cdots \Mod{u_{i-1}'}\,\Mod{u_i'-u_i}\, \Mod{u_{i+1}}\cdots \Mod{u_k}\;. \]
\end{prop}

\dem
Take the sequence of multi-vectors defined by
$u^{(i)}=( u_1',\ldots, u_{i}', u_{i+1},\ldots, u_k)\in G^k$\,
($0\leq i \leq k$). Notice that $u^{(0)}=u$ and $u^{(k)}=u'$.
The stated inequality follows from
\begin{align*}
\Delta^k_{u'} f(x) - \Delta^k_{u} f(x) &= 
\sum_{i=1}^k ( \Delta^k_{u^{(i)}}  - \Delta^k_{u^{(i-1)}})\, f(x)\\
&= 
\sum_{i=1}^k ( \Delta_{u_i'}-\Delta_{u_i} )\,
 \Delta^{k-1}_{(u_1',\ldots, u_{i-1}',u_{i+1},\ldots, u_k) }\, f(x)  \\
&= 
\sum_{i=1}^k   \Delta_{u_i'-u_i} \, \tau_{u_i}\, 
 \Delta^{k-1}_{(u_1',\ldots, u_{i-1}',u_{i+1},\ldots, u_k) }\, f(x)  \\
&= 
\sum_{i=1}^k   \Delta_{u_i'-u_i} \,
 \Delta^{k-1}_{(u_1',\ldots, u_{i-1}',u_{i+1},\ldots, u_k) }\, f(x+u_i)  \;.
\end{align*}
\cqd

\bigskip

\section{Open Domains}

Let $\Omega$ be an open  convex set in $X$.
We denote by $\CD{k-1, 1}{\Omega}{Y}$ the space of class $C^{k-1}$
functions $f:\Omega\to Y$ whose $(k-1)$-derivative $D^{k-1}f:\Omega\to \Linsim{k-1}{X}{Y}$,
$x\mapsto D^{k-1} f_x$
is a Lipschitz function. We denote by ${\rm Lip}(D^{k-1}f)$ the Lipschitz constant
of this function.

\begin{teor}\label{teor:B}
\; 
$L^k(\Omega,Y)=\CD{k-1, 1}{\Omega}{Y}$ and 
\, $\spnrm{ \Delta^k f}_\Omega = {\rm Lip}(D^{k-1}f)$
for every function $f\in Y^\Omega$
in this space.
\end{teor}

\bigskip

To prove this theorem we need the following

\begin{prop}
\label{prim}
Assume $\Omega=D(x_0,r)$ is a disk with center $x_0$ and radius $r$.
Given $f\in L^{k+1}(\Omega,Y)$\, there are \, $h\in\Pol{k}{X}{Y}$\,
and \, $g\in\CD{k}{\Omega}{Y}$ such that\, $f=h+g$\, and 
\begin{enumerate}
\item[i)]\,\, $D^i g(x_0)=0$, for each $i=0,1,\ldots, k$,
\item[ii)]\,  $\nrm{ D^i g}\leq r^{k+1-i}\,\spnrm{ \Delta^k f}_\Omega $ 
over $\Omega$, for each $i=0,1,\ldots, k$,
\item[iii)]  ${\rm Lip}(D^{k}g)\leq \spnrm{ \Delta^k f}_\Omega $.
\end{enumerate}
\end{prop}

\bigskip

This proposition is proved at the end of this section.
The polynomial $h$ above is unique.
It is the $k^{\text{th}}$ Taylor polynomial of $f$ at  $x_0$.

\bigskip

\demof{ Proof of theorem~\ref{teor:B} }
Let us prove $L^{k+1}(\Omega,Y)\subseteq \CD{k, 1}{\Omega}{Y}$.
Apply proposition~\ref{prim} to every disk $D(x_0,r)$
contained in $\Omega$. We conclude that $f$ is of class $C^{k}$
and $D^{k} f$ has Lipschitz constant $\leq \spnrm{\Delta^k f}_\Omega$ on every disk
$D(x,r)\subseteq \Omega$.
It follows  that $D^{k} f$ has Lipschitz constant $\leq \spnrm{\Delta^k f}_\Omega$  
over the domain $\Omega$. Take $x,y\in\Omega$.
We can decompose the line segment $[x,y]\subseteq \Omega$
in points $x=x_0, x_1,\ldots, x_n= y$ such that 
$\nrm{y-x}=\sum_{i=1}^{n} \nrm{x_{i}-x_{i-1}}$, and each pair
$x_{i-1},\,x_i$ is contained in a common disk $D(x,r)\subseteq \Omega$.
Then
\begin{eqnarray*}
\nrm{ D^{k-1}f(y)- D^{k-1}f(x) } &\leq& \sum_{i=1}^n
\nrm{ D^{k-1}f(x_i)- D^{k-1}f(x_{i-1}) }\\
&\leq & 
\sum_{i=1}^n \spnrm{\Delta^k f}_\Omega\,\nrm{x_{i}-x_{i-1}} 
= \spnrm{\Delta^k f}_\Omega\,\nrm{y-x}\;.
\end{eqnarray*}

Let us now prove $\CD{k, 1}{\Omega}{Y}\subseteq L^{k+1}(\Omega,Y)$.
Assume $f\in \CD{k, 1}{\Omega}{Y}$. 
Applying proposition~\ref{mean value theorem} to $\Delta^{k} f(x)$,
\begin{eqnarray*}
\nrm{\Delta_v\Delta^{k}_u f(x)} &=& 
\nrm{ \Delta^{k}_u f(x+v) - \Delta^{k}_u f(x)}\\
&\leq & 
\int_{[0,1]^{k}} \nrm{ D^{k}_u f(x + v +t\cdot u) - D^{k}_u f(x+t\cdot u )}\,dt \\
&\leq & {\rm Lip}(D^k f)\,\nrm{u}^{k}\,\nrm{v}\;,
\end{eqnarray*}
which proves that $f\in L^{k+1}(\Omega,Y)$ with $\spnrm{\Delta^{k+1} f}_\Omega\leq {\rm Lip}(D^k f)$.
The convexity of $\Omega$ is needed by the assumption of proposition~\ref{mean value theorem}.
Notice that the norm in $\Linsim{k}{\R^d}{Y}$ is such that
for $\xi,\eta \in\Linsim{k}{\R^d}{Y}$, and $u\in X^{k}$,\;
$\nrm{\xi(u)-\eta(u)}\leq \nrm{\xi-\eta}\,\nrm{u}^{k}$.
\cqd

\bigskip

Let $f\in L^{k+1}(\Omega,Y)$.
We define a $k$-algebraic derivative at each point $x\in \Omega$,\;
$\nabla^k f(x)\in\Linsim{k}{\R^d}{Y}$.
These derivatives are obtained recursively but top-bottom, instead of the
usual bottom-top infinitesimal approach.
We define this derivative as a net limit of finite differences of order
$k$ of the function $f$.

\bigskip

Given $n\in\NN^k$, $u\in X^k$, we set
$\frac{u}{n}=\left( \frac{u_1}{n_1},\ldots, \frac{u_k}{n_k}\right)\in X^k$. 
Given $j\in [n]$, we write
$\frac{j}{n}=\left( \frac{j_1}{n_1},\ldots, \frac{j_k}{n_k}\right)\in \Q^k$
and $\Mod{\frac{j}{n}} = \frac{j_1}{n_1} + \ldots + \frac{j_k}{n_k}$.
A simple computation shows that
\begin{lema}
\label{smjn}
For every $k\in\N$ and $n\in\N^k$ with $n_i>0$ \; $(1\leq i\leq k)$,
\begin{equation}
\label{sumjn}
 \sum_{j\in [n]} \Mod{\frac{j}{n}}=\frac{\ovl{n}}{2}\,\left(k-\frac{1}{n_1}-\ldots -\frac{1}{n_k} \right) 
\; \leq \; \frac{\ovl{n}\,k}{2}\;.
\end{equation}
\end{lema}

\bigskip

\begin{lema}
\label{fdkr}
Let $f\in L^{k+1}(\Omega,Y)$.
Then for every $n\in \N^k$,
$$ \nrm{ \Delta^k_u f(x) - \ovl{n}\,\Delta^k_{\frac{u}{n}} f(x) }\leq \frac{k}{2}\,
\spnrm{\Delta^{k+1} f}_\Omega\,
(\max_i \nrm{u_i})\,\nrm{u}^k\;.$$
\end{lema}

\dem
Using~(\ref{monot}) and~(\ref{sumjn}), we have, for a given $x\in\Omega^{u}$,
\begin{eqnarray*}
\nrm{ \Delta^k_{u} f(x) - \ovl{n}\,\Delta^k_{\frac{u}{n}} f(x) } &\leq &
\sum_{j\in [n] }
\nrm{ \Delta^k_{\frac{u}{n}} f\left(x+\frac{j}{n}\cdot u\right) - \Delta^k_{\frac{u}{n}} f(x) } \\
&\leq & \sum_{j\in [n] }
\nrm{ \Delta_{\frac{j}{n}\cdot u} \Delta^k_{\frac{u}{n}} f(x) } \\
&\leq & \sum_{j\in [n] }
\spnrm{\Delta^{k+1} f}_\Omega \,\nrm{ \frac{j}{n}\cdot u}\,\nrm{\frac{u}{n}}^k\\
&\leq & \spnrm{\Delta^{k+1} f}_\Omega\,\frac{\nrm{u}^k}{\ovl{n}}  \, \sum_{j\in [n] }
\Mod{ \frac{j}{n} }\,(\max_i \nrm{u_i})\\
&\leq & \frac{k}{2}\,\spnrm{\Delta^{k+1} f}_\Omega\,(\max_i \nrm{u_i}) \,\nrm{u}^k\;.
\end{eqnarray*}
\cqd

\begin{lema}
\label{des_cauchy}
Let $f\in L^{k+1}(\Omega,Y)$.
Then for every $n,\, p\in \N^k$,
$$ \nrm{ \ovl{p}\,\Delta^k_{\frac{u}{p}} f(x) - \ovl{p}\,\ovl{n}\,\Delta^k_{\frac{u}{n\,p}} f(x) }
\leq \frac{k}{2}\,\spnrm{\Delta^{k+1} f}_\Omega\,
\left(\max_i \nrm{\frac{u_i}{p_i}} \right)\,\nrm{u}^k\;.$$
\end{lema} 

\dem
This follows by lemma~\ref{fdkr},
replacing $u$ by $\frac{u}{p}$ and multiplying both sides by $\ovl{p}$.
\cqd

\bigskip

Let $(\N^k,\succeq)$ be the partial ordered set of natural numbers with the order\; $\succeq$
$$m\succeq n \quad \Leftrightarrow \quad 
m = n\ast s  \; \text{ for some } s\in\N^k \;,$$
where $n\ast s=(n_1 s_1,\ldots, n_k s_k)$ for $n=(n_1,\ldots, n_k)$
and $s=(s_1,\ldots, s_k)$.

We denote by $\slim{n} (x_n)$,
when it exists, the limit of a convergent net $(x_n)_{n\in (\N^k,\succ)}$.

\bigskip

\begin{prop}
\label{cauchy}
Given $f\in L^{k+1}(\Omega,Y)$, 
 $u\in X^{k}$ and $x\in\Omega^{u}$,\,
then \; $\left(\;\ovl{n}\,\Delta^k_{\frac{u}{n}} f(x)\; \right)_{n\in (\N^k,\succ)}$\; 
is a convergent net.
\end{prop}

\dem
From lemma~\ref{des_cauchy} above, this net is Cauchy,
and since $Y$ is Banach, it must converge.
\cqd

\bigskip

Let $f\in L^{k+1}(\Omega, Y)$.
By proposition~\ref{cauchy},
we can define, for any $u\in X^k$ and  $x\in \Omega$,
\begin{equation}
\label{gradk}
\nabla^k_u f(x) := \, \slim{n}  \ovl{n}\,\Delta^k_{ \frac{u}{n} } f(x)\;.
\end{equation}  
Notice that since $\Omega$ is open, $x\in \Omega^{\frac{u}{n}}$ for all
sufficiently large $n$.
We also define 
$$\spnrm{\nabla^k_u f}_\Omega := \sup_{x\in \Omega} \nrm{\nabla^k_uf(x)}
\;\in [0,+\infty] \;\;\text{ and } $$
$$\spnrm{\nabla^k f}_\Omega := \sup_{u\in (X-\{0\})^k} 
\frac{ \nrm{\nabla^k_uf(x)} }{\nrm{u}^k}
\;\in [0,+\infty]\;.$$
\bigskip

\begin{prop}
\label{nablaprop}
Given  $f,g\in L^{k+1}(\Omega,Y)$\, and\,  $x\in\Omega$,
\begin{enumerate}
\item[(a)] the map $\nabla^k f(x):X^k\rightarrow Y$,\;
$u\mapsto \nabla^k_u f(x)$,\, belongs to $\Linsim{k}{\R^d}{Y}$\;.
\item[(b)] $\spnrm{\Delta^k f}_\Omega\, = \, \spnrm{\nabla^k f}_\Omega$.

\item[(c)]
$\displaystyle{ \Delta^k_u (f+g)(x)= \Delta^k_u f(x)+\Delta^k_u g(x) }$
\; and \;
$\displaystyle{ \nabla^k_u (f+g)(x)= \nabla^k_u f(x)+\nabla^k_u g(x) }.$
\end{enumerate}
\end{prop}

\dem
\begin{enumerate}
\item[(a)] The symmetry is a consequence of the commutativity of the difference
operators $\Delta_u$. Therefore, it is enough to show the first argument additivity.
Given $u,v\in X$, and $w\in X^{k-1}$,
\begin{eqnarray*}
\nabla^k_{(u+v),w} f &=&
\slim{n,p} p\,\ovl{n}\, \Delta_{\frac{u+v}{p}}\Delta^{k-1}_{\frac{w}{n}} f\\
&=&
\slim{n,p} p\,\ovl{n}\,\left[
\Delta_{\frac{u}{p}}\Delta^{k-1}_{\frac{w}{n}} f
+
\tau_{\frac{u}{p}}\,\Delta_{\frac{v}{p}}\Delta^{k-1}_{\frac{w}{n}} f
\right]\\
&=&
\slim{n,p} p\,\ovl{n}\,\left[
\Delta_{\frac{u}{p}}\Delta^{k-1}_{\frac{w}{n}} f
+
\Delta_{\frac{v}{p}}\Delta^{k-1}_{\frac{w}{n}} f
-
\Delta_{\frac{u}{p}}\Delta_{\frac{v}{p}}\Delta^{k-1}_{\frac{w}{n}} f
\right]\\
&=&
\nabla^k_{u,w} f + \nabla^k_{v,w} f\;.
\end{eqnarray*}
We note that, since $f\in L^{k+1}(\Omega,Y)$,
$$ p\,\ovl{n}\,\nrm{\; \Delta_{\frac{u}{p}}\Delta_{\frac{v}{p}}\Delta^{k-1}_{\frac{w}{n}} f\; } \leq
\frac{\spnrm{\Delta^{k+1} f}_\Omega}{p}\, \nrm{u}\,\nrm{v}\,\nrm{w}^{k-1}\;,$$
which proves that
$$\slim{n,p} p\,\ovl{n}\,\Delta_{\frac{u}{p}}\Delta_{\frac{v}{p}}\Delta^{k-1}_{\frac{w}{n}} f =0\;.$$

\item[(b)] The inequality $ \spnrm{\nabla^k f}_\Omega\leq
\spnrm{\Delta^k f}_\Omega$ holds because
$$ \frac{ \nrm{ \ovl{n}\,\Delta^k_{\frac{u}{n}} f(x) } }{\nrm{u}^k }
= \frac{ \nrm{ \Delta^k_{\frac{u}{n}} f(x) } }{\nrm{\frac{u}{n}}^k }
\leq \spnrm{\Delta^k f}_\Omega\;,
$$
for every $u\in X^k$, \, $x\in \Omega^u$,\, and\, $n\in\N^k$.
The reverse inequality follows from proposition~\ref{monot:prop2} bellow.

\item[(c)] This is a direct consequence of the linearity of the difference operators
$\Delta_u$.
\end{enumerate}
\cqd

\begin{prop}
\label{monot:prop2}
For every  $f\in Y^\Omega$ and $u\in X^{k}$,\, the net
$$\left(\;\ovl{n}\,\spnrm{ \Delta^k_{\frac{u}{n}} f}_\Omega\; \right)_{n\in (\N^k,\succ)}$$
is monotonous increasing, i.e.,\;
$n\succeq m$ \quad $\Rightarrow$ \quad
$\ovl{n}\,\spnrm{ \Delta^k_{\frac{u}{n}} f}_\Omega \; \geq \;
\ovl{m}\,\spnrm{ \Delta^k_{\frac{u}{m}} f}_\Omega$ \;.
\end{prop}

\dem
Follows from proposition~\ref{monot:prop1}.
\cqd

\bigskip

\begin{prop}
\label{reducao}
Let $\Omega=D(x_0,r)$ be a disk with center $x_0$ and radius $r$.
Given $f\in L^{k+1}(\Omega,Y)$
the function $g\in Y^\Omega$, \;
$$g(x)= f(x)- \frac{1}{k!}\diag (\nabla^k f(x_0) )(x-x_0)\;,$$
is in $L^k(\Omega,Y)$\, with \, $\spnrm{\Delta^k g}_\Omega \leq r\, \spnrm{\Delta^{k+1} f}_\Omega$\, and \,
$\nabla^k g(x_0)=0$.
\end{prop}

\dem
We have $g\in L^{k+1}(\Omega,Y)$ 
because $f\in L^{k+1}(\Omega,Y)$
and  
$$p(x)=\frac{1}{k!}\diag (\nabla^k f(x_0) )(x-x_0) $$
is an homogeneous polynomial of degree $k$, and so,
by proposition~\ref{derivn}(2),
$\Delta^{k+1}_u p\equiv 0$.
It is also clear that $\spnrm{\Delta^{k+1} g}_\Omega = \spnrm{\Delta^{k+1} f}_\Omega$.
By proposition~\ref{derivn}(1),
$\Delta^k_u p(x)\equiv \nabla^k_u f(x_0)$, which implies that
$\nabla^k_u p(x)\equiv \nabla^k_u f(x_0)$. 
Therefore, by proposition~\ref{nablaprop}(c),
$\nabla^k g(x_0)=0$. Because $\spnrm{\Delta^{k+1} g}_\Omega = \spnrm{\Delta^{k+1} f}_\Omega$,
we have, for every $u\in X$, $w\in X^k$ and $x\in\Omega^{u,w}$,
$$
\ovl{n}\, \nrm{ \Delta_u \Delta^k_{\frac{w}{n}} g(x) }\leq
\spnrm{\Delta^{k+1} f}_\Omega\,\nrm{u}\,\nrm{w}^k \;. $$
Taking the limit when $n\rightarrow\infty$,
we obtain
$$
 \nrm{ \Delta_u \nabla^k_{w} g(x) } \leq
\spnrm{\Delta^{k+1} f}_\Omega\,\nrm{u}\,\nrm{w}^k \;.$$
Replacing $x$ by $x_0$, and $u$ by $x-x_0$,
we have
$$
 \nrm{ \nabla^k_{w} g(x) } \leq
\spnrm{\Delta^{k+1} f}_\Omega\,\nrm{x-x_0}\,\nrm{w}^k
\leq
r\,\spnrm{\Delta^{k+1} f}_\Omega\,\nrm{w}^k\;,$$
for all $w\in X^k$ and $x\in\Omega^{w}$.
Therefore, by proposition~\ref{nablaprop}(b),
$$\spnrm{\Delta^k g}_\Omega =\spnrm{\nabla^k g}_\Omega\leq  r\,\spnrm{\Delta^{k+1} f}_\Omega\;,$$
which in particular proves that $g\in L^k(\Omega,Y)$.
\cqd

\bigskip
\bigskip

\demof{Proof of proposition~\ref{prim}}
Consider the following claim, which depends on the parameter $s\in\N$:

\smallskip

$\mathcal{S}(s)$\quad $:\Leftrightarrow$\quad
there exists $h_s\in \oplus_{i=s}^k \Hom{i}{X}{Y}$ 
such that $g_s = f- h_s$ satisfies:
\quad\begin{enumerate}
\item[i)]\,\, $\nabla^i g_s(x_0)=0$, for each $i=s,\ldots, k$,
\item[ii)]\,  $\spnrm{ \nabla^i g_s}_\Omega\leq r^{k+1-i}\,\spnrm{\Delta^{k+1} f}_\Omega$, 
for each $i=s,\ldots, k$,
\item[iii)] $\spnrm{ \Delta_v \nabla^k_u g_s }_\Omega \leq 
\spnrm{\Delta^{k+1} f}_\Omega\,\nrm{v}\,\nrm{u}^k$ 
\end{enumerate}
Note that all the derivatives above are defined in the weak sense~(\ref{gradk}).
We shall prove, by regressive induction, the claim
$\mathcal{S}(0)$.
By proposition~\ref{reducao}, the claim
$\mathcal{S}(k)$ holds.
It lefts to show that \; 
$$ \mathcal{S}(s)\; \Rightarrow \mathcal{S}(s-1),\quad
\text{for any }\; 1\leq s\leq k\;.$$
Assume that 
$\mathcal{S}(s)$ holds and let us apply again
proposition~\ref{reducao} to the function $g_s$. 
By $\mathcal{S}(s)$, we have
$\spnrm{\Delta^{s} g_s}_\Omega\leq r^{k+1-s} \spnrm{\Delta^{k+1} f}_\Omega$.
Therefore, proposition~\ref{reducao} guarantees the existence of a homogeneous polynomial
~$p_{s-1}\in\Hom{s-1}{X}{Y}$
such that
\begin{enumerate}
\item[a)]
$\spnrm{\Delta^{s-1} g_{s-1}}_\Omega\leq r^{k+2-s} \spnrm{\Delta^{k+1} f}_\Omega$,
where $g_{s-1}=g_s-p_{s-1}$, \, and
\item[b)] $\nabla^{s-1} g_{s-1}(x_0)=0$\;.
\end{enumerate}
We define\; $h_{s-1}=h_s+p_{s-1}$\;.
Then
\begin{eqnarray*}
g_{s-1} &=& g_s-p_{s-1} = (f-h_s)- p_{s-1}\\
&=& f-(h_s+p_{s-1}) = f- h_{s-1}\;.
\end{eqnarray*}
By proposition~\ref{derivn},
$\nabla^i p_{s-1}\equiv 0$, for all $i\geq s$.
Therefore, considering the induction hypothesis, $\mathcal{S}(s)$,
we have for each $i\geq s$,
\quad\begin{enumerate}
\item[i)]\,\, $\nabla^i g_{s-1}(x_0)=\nabla^i g_s(x_0) +
\nabla^i p_{s-1}(x_0) = 0$, 
\item[ii)]\,  $ \spnrm{ \nabla^i g_{s-1} }_\Omega
= \spnrm{ \nabla^i g_{s} }_\Omega
\leq  r^{k+1-i}\, \spnrm{\Delta^{k+1} f}_\Omega$,
\item[iii)] $ \spnrm{ \Delta_v \nabla^k_u g_{s-1} }_\Omega 
=\spnrm{ \Delta_v \nabla^k_u g_{s} }_\Omega 
\leq \spnrm{\Delta^{k+1} f}_\Omega\,\nrm{v}\,\nrm{u}^k$. 
\end{enumerate}
Finally,
a), b), i), ii) e iii) show that
$\mathcal{S}(s-1)$ holds.
\cqd

\bigskip
\bigskip

\section{Lattice Domains}

In this section we assume $X=\R^d$ and
$\Gamma=m_1^{-1}\Z\times\ldots \times m_d^{-1}\Z\subset\R^d$ for some $m\in\Z^d_+$.
We write $e_i=(0,\ldots, 0,m_i^{-1},0,\ldots, 0)$ ($1\leq i\leq d$),
so that  $\Escr=\{ e_1,\ldots, e_d \}$ is a basis of both $\Gamma$ (as free abelian group) and  $\R^d$.
We fix the sum-norm \, $\nrm{(x_1,\ldots, x_d)}:=\sum_{i=1}^d \Mod{x_i}$\, in $\R^d$.
Notice that for every $(c_1,\ldots, c_d)\in\R^d$ we have
\begin{equation}\label{sum:norm}
 \nrm{\sum_{i=1}^d c_i\,e_i} = \sum_{i=1}^d \Mod{c_i} \,\nrm{e_i}\;.
\end{equation}
Theorem A is about extending lattice functions $f:\Gamma/\Z^d\to Y$ with bounded $k$-differences
to functions in $\CD{k-1,1}{\T^d}{Y}$ with Lipschitz $(k-1)$-derivatives.
This amounts to extend a $\Z^d$-periodic function  $f:\Gamma \to Y$ to a 
$\Z^d$-periodic function $F:\R^d\to Y$. We say that a function $f:\Gamma\to Y$ is $\Z^d$-periodic \, iff\,
$f(x)=f(x+h)$ for all $x\in\Gamma$  and $h\in \Z^d$, and denote by $\Fscr_{\Z^d}(\Gamma,Y)$ the set of all  
$\Z^d$-periodic functions $f:\Gamma\to Y$.
In this context the seminorms $\spnrm{\Delta^k f}_\Gamma$ are always well-defined
because $\Gamma/\Z^d$ is a finite set (a group actually) with $m_1\cdots m_d$ elements.
We set
$$ L^k_{\Z^d}(\Gamma,Y):=\{\, f\in \Fscr_{\Z^d}(\Gamma,Y)\,:\,
\spnrm{\Delta^k f}_\Gamma<+\infty \, \}\;.$$
Vectors in $B = \{ e_1,\ldots, e_d, -e_1,\ldots, -e_d\}$ will be referred as 
{\em basic vectors}. We define $B^k = \{ e_1,\ldots, e_d, -e_1,\ldots, -e_d\}^k$.
Elements of $B^k$ will be referred as {\em basic multi-vectors}.
Given $f\in\Fscr_{\Z^d}(\Gamma, Y)$ and any subset $\Sigma\subset \Gamma$ we set
\begin{align*}
 \spnrm{\Delta^k_u f}_\Sigma &:= 
\sup_{x\in \Sigma } \nrm{ \Delta^k_u f(x)}
\;\in [0,+\infty]\,, \; \text{ and }\\
 \spnrm{\Delta^k f}_\Sigma &:= 
\sup_{ u\in (\Gamma-\{0\})^k} \frac{ \spnrm{\Delta^k_{u} f}_\Sigma }{\nrm{u}^k}
\;\in [0,+\infty] \;.
\end{align*}
We can restrict the computation of difference norms to basic multi-vectors.

\begin{prop} \label{nrm:char} Given  $f\in \Fscr_{\Z^d}(\Gamma,Y)$
and $\Sigma\subset \Gamma$,
$$ \spnrm{\Delta^k f}_\Sigma = 
\max_{ u\in B^k}  \frac{ \spnrm{\Delta^k_{u} f}_\Sigma }{\nrm{u}^k}  \;.$$
\end{prop}

\dem
The inequality $\geq $ is obvious. 
For the reverse inequality we use proposition~\ref{gen:monot} below  and~(\ref{sum:norm}).
\cqd


\begin{prop}[General Monotonicity]
\label{gen:monot}
Given $f\in \Fscr_{\Z^d}(\Gamma, Y)$, 
$\Sigma\subset \Gamma$,
$\,u \in \Gamma^d\,$ and 
a $k\times d$ matrix $\,(n_{i,j}) \,$ with entries in $\NN$,
\begin{equation*}
\spnrm{ \Delta^k_{\left( 
\sum_{j_1=1}^d n_{1,j_1} u_{j_1},\, \cdots,\,
\sum_{j_k=1}^d n_{k,j_k} u_{j_k}  
\right) } f }_\Sigma  \leq
\sum_{j_1=1}^d
\cdots
\sum_{j_k=1}^d n_{1,j_1}\cdots n_{k,j_k} \spnrm{ \Delta^k_{(u_{j_1},\cdots, u_{j_k})} f}_\Sigma 
 \;.
\end{equation*}
\end{prop}

\bigskip

\dem
It follows from the definition of the difference
operators $\Delta_u$ that
$$\Delta_{\sum_{j=1}^d u_j}  =
\sum_{j=1}^d \tau_{u_1+ \cdots + u_{j-1}}\circ \Delta_{u_j} \;. $$
Therefore,
composing
$\left( \Delta_{\sum_{j_1=1}^d u_{1,j_1}}\right) \circ \cdots \circ 
\left( \Delta_{\sum_{j_k=1}^d u_{k,j_k}}\right)  $
we obtain
$$\spnrm{ \Delta^k_{\left( 
\sum_{j_1=1}^d u_{1,j_1},\, \cdots,\,
\sum_{j_k=1}^d u_{k,j_k}  
\right) } f }_\Sigma  \leq
\sum_{j_1=1}^d
\cdots
\sum_{j_k=1}^d \spnrm{ \Delta^k_{(u_{1,j_1},\cdots, u_{k,j_k})} f}_\Sigma
\;.$$

The inequality follows replacing $u_{i,j_i}$ by $n_{i,j_i}\,u_{j_i}$, 
and using~(\ref{monot}).
\qed

\bigskip

We consider on the space $\Linsim{k}{\R^d}{Y}$ of (continuous) symmetric $k$-multi-linear functions
the usual norm
\[ \nrm{\xi} = \sup_{u\in (\R^d-\{0\})^k} \frac{ \nrm{\xi(u)} }{\nrm{u}^k } \;.\]
By~(\ref{sum:norm}) this norm is attained at basic multi-vectors,\,
$\nrm{\xi} = \max_{ u\in B^k}    { \nrm{\xi(u)} }/{\nrm{u}^k } $.

\bigskip

We now associate to each function $f\in L^k_{\Z^d}(\Gamma, Y)$ an approximate 
$k^{\text{th}}$ order derivative at each point $x\in  \Gamma$
 by averaging  differences of order $k$ around $x$.
Given $x\in  \Gamma$ and $\alpha \in \{0,1\}^d$, we define
$\Theta^k_{\alpha}(f)(x) \in \Linsim{k}{\R^d}{Y}$ by
\begin{equation}
\label{thetak}
\Theta^k_{\alpha}(f)(x)((-1)^{\alpha_{i_1}} e_{i_1}, \ldots, (-1)^{\alpha_{i_k}} e_{i_k})
:=\Delta^k_{(\,(-1)^{\alpha_{i_1}} e_{i_1},\, \ldots\, , \, (-1)^{\alpha_{i_k}} e_{i_k}\,)}f(x),
\end{equation}
for each multi-index $\iota=(i_1,\ldots, i_k)\in \{1,\ldots, d\}^k$.
The $k$-multi-linear map $\Theta^k_{\alpha}(f)(x)$ is well-defined because
$\{(-1)^{\alpha_1}e_1, \cdots, (-1)^{\alpha_d}e_d\}$ is a basis for $\R^d$.
Then we define $\Theta^k(f)(x) \in \Linsim{k}{\R^d}{Y}$ by averaging
$$\Theta^k(f)(x):=\frac{1}{2^d} \sum_{\alpha \in \{0,1\}^d} \Theta^k_{\alpha}(f)(x)\;.$$
Notice that for $k=0$ and $\alpha\in\{0,1\}^d$, \, 
$\Theta^0(f)(x)=\Theta^k_{\alpha}(f)(x) = f(x)$.

We also define for each $u\in \Gamma^k$ and each $\Sigma\subset \Gamma$,
$$\spnrm{\Theta^k_u f }_\Sigma := \sup_{x \in  \Sigma} \nrm{\Theta^k f(x) (u) }
\;\in [0,+\infty] \;\;\text{ and } $$
$$\spnrm{\Theta^k f}_\Sigma := \sup_{u\in (\Gamma-\{0\})^k}
\frac{ \spnrm{\Theta^k_u f}_\Sigma}{\nrm{u}^k}
\;\in [0,+\infty]\;.$$
\bigskip

\bigskip

\begin{prop}
\label{Thetaprop}
Given
 $\,f,g\in \,L^{k}_{\Z^d}(\Gamma,Y)$ \, and \,  $x\in\Gamma$,
\begin{enumerate}
\item[(a)] the map $\Theta^k f(x):(\R^d)^k\rightarrow Y$,\;
$u\mapsto \Theta^k f(x)(u)$,\, belongs to $\Linsim{k}{\R^d}{Y}$\;.

\item[(b)]
$\displaystyle{ \Theta^k_u (f+g)(x)= \Theta^k_u f(x)+\Theta^k_u g(x) }.$
\end{enumerate}
\end{prop}

\bigskip

\dem
\begin{enumerate}
\item[(a)] The symmetry is a consequence of the commutativity of the difference
operators $\Delta_u$.

\item[(b)] This is a direct consequence of the linearity of the difference operators
$\Delta_u$.
\end{enumerate}
\cqd

\bigskip

\begin{prop} 
\label{diff:sim}
Given  $\,f\in Y^\Gamma$, $u\in (\R^d)^k$ and $\alpha\in \{0,1\}^k$,
$$ \Delta^k_{(-1)^\alpha\ast u} f (x) = (-1)^{\Mod{\alpha}}\,\Delta^k_u f (x-\alpha\cdot u)
\;,$$
where $(-1)^\alpha\ast u = (\, (-1)^{\alpha_1} u_1,\,\ldots,\, (-1)^{\alpha_k} u_k\,)$.
\end{prop}

\dem
The proof goes by induction in $k$.
\cqd

\bigskip

\begin{prop} 
\label{Theta:formula}
Given  $\,f\in Y^\Gamma$, $x\in  \Gamma$ and $(i_1,\ldots, i_k)\in\{1,\ldots, d\}^k$,
$$ \Theta^k(f)(x)(e_{i_1},\ldots, e_{i_k}) 
=\frac{1}{2^k}\, \Delta^k_{(\,2 e_{i_1},\,\ldots, \,2 e_{i_k}\,)} f
\left(\,x-e_{i_1} -\ldots - e_{i_k}\,  \right) $$
\end{prop}

\dem
The first two steps follow from the definitions.
The third step uses proposition~\ref{diff:sim}.
The sum obtained can be grouped in $2^k$ groups with
$2^{d-k}$ equal summands each, which justifies the fourth and final step.

\begin{align*}
\Theta^k(f)(x)(e_{i_1},\ldots, e_{i_k}) &=
\frac{1}{2^d} \sum_{\alpha \in \{0,1\}^d} \Theta^k_{\alpha}(f)(x)(e_{i_1},\ldots, e_{i_k})\\
 &=
\frac{1}{2^d} \sum_{\alpha \in \{0,1\}^d} (-1)^{\alpha_{i_1}+\cdots + \alpha_{i_k}}
\Delta^k_{(\,(-1)^{\alpha_{i_1}} e_{i_1},\, \ldots\, , \, (-1)^{\alpha_{i_k}} e_{i_k}\,)}f(x)\\
&=
\frac{1}{2^d} \sum_{\alpha \in \{0,1\}^d}  
\Delta^k_{(\, e_{i_1},\, \ldots\, , \,  e_{i_k}\,)}
f\left(\, x -\alpha_{i_1} e_{i_1}-\cdots - \alpha_{i_k} e_{i_k}\,\right)\\
&=
\frac{1}{2^k} \sum_{\alpha \in \{0,1\}^k}  
\Delta^k_{(\, e_{i_1},\, \ldots\, , \,  e_{i_k}\,)}
f\left(\, x -\alpha_{1} e_{i_1}-\cdots - \alpha_{k} e_{i_k}\,\right)\\
\end{align*}
\cqd

\bigskip

\begin{prop}\label{Theta:polynomial}
Given $f \in \Pol{m}{\R^d}{Y}$, $x,x_0\in X$, $u=(u_1,\ldots, u_k)\in (\R^d)^k$  and $0\leq k \leq m$,
we have
\begin{align*}
\Theta^k_u f(x)- D^k_u f(x)
 &=  \sum_{r_0+r_1+\cdots + r_k=\tilde{r}}
 (-1)^{\tilde{r}-k +r_0}\,\frac{\tilde{r}!}{r_0!\, r_1!\,\ldots\, r_k!} 
D^{\tilde{r}}_{( (x-x_0)^{(r_0)}, u_1^{(r_1)},\ldots, u_k^{(r_k)} ) } f(x_0)
\end{align*}
where the sum is taken over all $(r_0,r_1,\ldots, r_k)\in\N^{k+1}$ such that
 $r_i$ is odd for each $1\leq i \leq k$,  $r_i>1$ for some $1\leq i \leq k$
and $k+2\leq \tilde{r} \leq m$ where
 $\tilde{r}:= r_0+r_1+\cdots + r_k$.
\end{prop} 

\dem
This follows applying the following proposition to each term of the Taylor development of $f(x)$
at $x=x_0$.
\cqd

\begin{prop}\label{Theta:homogeneous:polynomial}
Given $f=\diag (\xi)\in \Hom{m}{\R^d}{Y}$, $u=(u_1,\ldots, u_k)\in (\R^d)^k$  and $0\leq k \leq m$,
we have
\begin{align*}
D^k_u f(x) &= \frac{m!}{(m-k)!} \xi( x^{(m-k)}, u_1,\ldots, u_k )\\
\Theta^k_u(f)(x) &=
D^k_u f(x) + \sum_{r_0+r_1+\cdots + r_k=m}
 (-1)^{m-k +r_0}\,\frac{m!}{r_0!\, r_1!\,\ldots\, r_k!} \xi( x^{(r_0)}, u_1^{(r_1)},\ldots, u_k^{(r_k)} ) 
\end{align*}
where the second sum is taken over all $(r_0,r_1,\ldots, r_k)\in\N^{k+1}$ such that
 $r_i$ is odd for each $1\leq i \leq k$,  $r_i>1$ for some $1\leq i \leq k$
and $r_0+r_1+\cdots + r_k=m$.
\end{prop}

\dem
For each $u\in (\R^d)^k$ we set \, $\tilde{u}=(x,u_1,\ldots, u_k)$,
which means $u_0=x$, and for $\alpha\in\{0,1\}^k$ we set\, $\tilde{\alpha}=(1,\alpha_1,\ldots, \alpha_k)$,
which means $\alpha_0=1$.
\begin{align*}
\Theta^k_u (f)(x) &= \frac{1}{2^k}\,
\sum_{\alpha\in\{0,1\}^k}
(-1)^{k-\Mod{\alpha}} f( x + 2\,\alpha\cdot u -\sum_{i=1}^k u_i ) \\
&= \frac{1}{2^k}\,
\sum_{\alpha\in\{0,1\}^k}
(-1)^{k-\Mod{\alpha}} f( x-\sum_{i=1}^k (-1)^{\alpha_i}\,u_i ) \\
&= \frac{1}{2^k}\,
\sum_{\alpha\in\{0,1\}^k}
(-1)^{k-\Mod{\alpha}} \xi( -\sum_{j_1=0}^k (-1)^{\alpha_{j_1}}\,u_{j_1},\ldots, 
-\sum_{j_m=0}^k (-1)^{\alpha_{j_m}}\,u_{j_m} ) \\
&= \frac{1}{2^k}\,
\sum_{\alpha\in\{0,1\}^k}
 \sum_{j_1=0}^k \ldots \sum_{j_m=0}^k 
(-1)^{k-\Mod{\alpha}+m+\alpha_{j_1}+\cdots + \alpha_{j_m}}
\xi(u_{j_1},\ldots,u_{j_m} ) \\
&= 
\sum_{r_0+r_1+\cdots + r_k=m}
c^{m,k}_{r_0,\ldots, r_k}\, \frac{m!}{r_0!\, r_1!\,\ldots\, r_k!} \xi( x^{(r_0)}, u_1^{(r_1)},\ldots, u_k^{(r_k)} ) 
\end{align*}
where for each $r=(r_0,r_1,\ldots, r_k)$ such that $r_0+r_1+\ldots + r_k=m$
\begin{align*} 
c^{m,k}_{r_0,\ldots, r_k} &:=
\frac{1}{2^k}\,\sum_{\alpha\in\{0,1\}^k}  (-1)^{k+m -\Mod{\alpha}+r\cdot\tilde{\alpha}}\\
&=
\frac{1}{2^k}\,\sum_{\alpha\in\{0,1\}^k}  (-1)^{m-k -\Mod{\alpha}+r_0 +\alpha_{1} r_1+\cdots + \alpha_{k} r_k} \\
&=
\frac{(-1)^{m-k +r_0}}{2^k}\, \sum_{\alpha\in\{0,1\}^k}  (-1)^{ \alpha_{1} (r_1-1)} \cdots (-1)^{ \alpha_{k} (r_k-1)} \\
&=
\frac{(-1)^{m-k +r_0}}{2^k}\, \left( \sum_{\alpha_1 \in\{0,1\}}  (-1)^{ \alpha_{1} (r_1-1)}\right)
\,  \cdots \, \left( \sum_{\alpha_k \in\{0,1\}}  (-1)^{ \alpha_{k} (r_k-1)}\right)  \\
&=
\frac{(-1)^{m-k +r_0}}{2^k}\, ((-1)^{r_1-1} +1) \,  \cdots \, ((-1)^{r_k-1} +1)  \\
&=\left\{\begin{array}{ccl}
(-1)^{m-k +r_0} &\text{ if } & r_i\; \text{ is odd for every }\; i \geq 1\\
0 &\text{ if } & r_i\; \text{ is even for some }\; i \geq 1\\
\end{array}\right.
\end{align*}
Finally, notice that $c^{m,k}_{r_0,\ldots, r_k}=1$ when $r_0=m-k$ and $r_i=1$ for every $i \geq 1$,
and the corresponding term is equal to $D^k_u f(x)=\frac{m!}{(m-k)!} \xi( x^{(m-k)}, u_1,\ldots, u_k )$.
\cqd

\bigskip

Given any subset $\Sigma\subset\Gamma$ define the neighbourhood of $\Sigma$,
$$N_k(\Sigma):=\{\, x+\sum_{i=1}^k u_i\,:\, x \in \Sigma, \,u=(u_1,\ldots, u_k)\in B^k\,\}\;.$$

\begin{prop}
\label{normatheta}
Given $f\in Y^\Gamma$ and $\Sigma\subset \Gamma$ finite, 
for every $\,x\in\Sigma$, $u \in \Gamma$ and  $\,w \in \Gamma^k$,
$$ \nrm{ \Delta_u \Theta^k_{w} f(x) } \leq
\spnrm{\Delta^{k+1} f}_{N_k(\Sigma)} \,\nrm{u}\,\nrm{w}^k \;.$$
\end{prop}

\bigskip

\dem
It is enough to consider the case where $u$ is a basic vector
and $w$ a basic multi-vector.
In this case the inequality follows from proposition~\ref{Theta:formula}
and inequality~(\ref{monot}).
\qed

\bigskip
\bigskip

We call {\em width } of the lattice $\Gamma$ to the positive number
$$\nrm{\Gamma}:= \max_{1\leq i\leq d} \nrm{e_i}=\left(\min_{1\leq i\leq d} m_i\right)^{-1}\;.$$

\begin{prop}
\label{Delta:Theta:ineq}
Given  $\,f\in Y^\Gamma$ and $\Sigma\subset \Gamma$ finite
\begin{enumerate}
\item[(a)] $\displaystyle  \spnrm{\Theta^k f}_\Sigma \, \leq \, \spnrm{\Delta^k f}_{N_k(\Sigma)} $,
\item[(b)] $\displaystyle  \spnrm{\Delta^k f}_{\Sigma}
\, \leq \, 
\spnrm{\Theta^k f}_\Sigma +  \frac{k}{2}\,\nrm{\Gamma}\,\spnrm{\Delta^{k+1} f}_{N_k(\Sigma)} $.
\end{enumerate}

\end{prop}

\dem
Item (a) follows from inequality~(\ref{monot}) and propositions~\ref{Theta:formula}
and~\ref{nrm:char}.
Let us now prove item (b).
Define $\xi_i=\{0,1\}^i\times\{0\}^{k-i}$ ($0\leq i \leq k$), so that
 $\xi_0\subset \xi_1\subset \ldots \subset \xi_k=\{0,1\}^k$.
Let us fix $\iota=(i_1,\ldots, i_k)\in \{1,\ldots, d\}^k$ and write
$e_\iota = (e_{i_1},\ldots, e_{i_k} )$.
By the calculation in the proof of proposition~\ref{Theta:formula}
we have 
$$ \Theta^k f(x)(e_\iota) = \frac{1}{2^k} \, 
\sum_{\alpha\in \xi_k} \Delta^k_{e_\iota}
f\left(\, x -\alpha\cdot e_\iota \,\right)  \;. $$
Define now for each $i=1,\ldots, k$,
\begin{align*}
 A_i &= \frac{1}{2^i} \,\left(
\sum_{\alpha\in \xi_{i-1}} \Delta^k_{e_\iota}
f\left(\, x -\alpha\cdot e_\iota \,\right)\, -\,
\sum_{\alpha\in \xi_{i}-\xi_{i-1}} \Delta^k_{e_\iota}
f\left(\, x -\alpha\cdot e_\iota\,\right)
\right) \\
&= - \frac{1}{2^i} \, 
\sum_{\alpha\in \xi_{i-1}} \Delta_{-e_i}\,\Delta^k_{e_\iota} f\left(\, x -\alpha\cdot e_\iota\,\right)
\;. 
\end{align*} 
Whence $\nrm{A_i}\leq \frac{1}{2}\,\spnrm{\Delta^{k+1} f}_\Gamma\,\nrm{e_i}\,\nrm{e_\iota}^k
\leq \frac{1}{2}\,\nrm{\Gamma}\,\spnrm{\Delta^{k+1} f}_\Gamma\,\nrm{e_\iota}^k $.
It is also easy to check that
$$\Delta^k_{e_\iota} f(x) = \Theta^k f(x)(e_\iota)  + A_k+\ldots + A_1\;,$$
and item (b) follows.
\cqd

\bigskip

Given $\Sigma\subset \Gamma$ and $x_0\in\Sigma$
define the {\em radius of } $\Sigma$ w.r.t. $x_0$ as
$$r_{x_0}(\Sigma):= \min\{\, \nrm{x-x_0}\,:\, x\in\Sigma\,\}\;.$$

\bigskip

\begin{prop}
\label{reducao2}
Given $f\in Y^\Gamma$, $\Sigma\subset \Gamma$  finite, and $x_0\in\Sigma$ 
consider the polynomial\,
$$p(x)=\frac{1}{k!}\diag \, (\Theta^k f(x_0) )(x-x_0)\,, $$
and the function $g \in Y^\Gamma$, \, $g(x)= f(x)- p(x)$.
Then\,
\begin{enumerate}
\item[(a)] \, $\Theta^k g(x_0)=0$  ,
\item[(b)] \, $\spnrm{\Delta^{k} g}_{\Sigma}
\leq  \left(\frac{k}{2}\,\nrm{\Gamma} +r_{x_0}(\Sigma) \right) \spnrm{\Delta^{k+1} f}_{N_k(\Sigma)}$ \, .  
\end{enumerate}
\end{prop}

\bigskip

\dem
Write $\Sigma'=N_k(\Sigma)$ and $r=r_{x_0}(\Sigma)$.
Clearly  $\spnrm{\Delta^{k+1}g}_{\Sigma'}=\spnrm{\Delta^{k+1}f}_{\Sigma'}$
because
$$p(x)=\frac{1}{k!}\diag (\Theta^k f(x_0) )(x-x_0)$$
is an homogeneous polynomial of degree $k$, and 
by proposition~\ref{derivn} (2) 
$\Delta^{k+1}_u p\equiv 0$.
By proposition~\ref{derivn} (1)
$\Delta^k_u p(x)\equiv \Theta^k_u f(x_0)$, which implies that
$\Theta^k_u p(x)\equiv \Theta^k_u f(x_0)$.
See also proposition~\ref{Theta:homogeneous:polynomial}.
Therefore, by proposition~\ref{Thetaprop} (b),
$\Theta^k g(x_0)=0$. It follows from proposition~\ref{normatheta} that
 for every $x\in\Sigma$, $u\in\Gamma$ and $w\in \Gamma^{k}$,
$$
 \nrm{ \Delta_u \Theta^k_{w} g(x) } \leq
\spnrm{\Delta^{k+1} f}_{\Sigma'} \,\nrm{u}\,\nrm{w}^k \;.$$
Replacing $x$ by $x_0$, and $u$ by $x-x_0$,
we have
$$
 \nrm{ \Theta^k_{w} g(x) } \leq
\spnrm{\Delta^{k+1} f}_{\Sigma'}\,\nrm{x-x_0}\,\nrm{w}^k
\leq
r\,\spnrm{\Delta^{k+1} f}_{\Sigma'} \,\nrm{w}^k,\;$$
for all $w \in \Gamma^{k}$ and $x \in \Sigma$.
This shows that $\spnrm{\Theta^k g}_{\Sigma}\leq r\,\spnrm{\Delta^{k+1} f}_{\Sigma'}$.
Therefore, by proposition~\ref{Delta:Theta:ineq},
\begin{align*}
\spnrm{\Delta^k g}_{\Sigma} & \leq
\spnrm{\Theta^k g}_{\Sigma} +\frac{k}{2}\, \spnrm{\Delta^{k+1} f}_{\Sigma'} \,\nrm{\Gamma}\\
& \leq r\,\spnrm{\Delta^{k+1} f}_{\Sigma'} +\frac{k}{2}\, \spnrm{\Delta^{k+1} f}_{\Sigma'}\,\nrm{\Gamma} \\
& = \left( r +\frac{k}{2}\,\nrm{\Gamma} \right)\, \spnrm{\Delta^{k+1} f}_{\Sigma'} \;.
\end{align*}
\cqd

\bigskip

In the next proposition we make use of the quantity
$$\ell(\Gamma):= \max_{1\leq i,j\leq d}  \frac{\nrm{e_i}}{\nrm{e_j}}
= \max_{1\leq i,j\leq d}  \frac{m_j}{m_i}  \;. $$

\begin{prop}\label{existence:px} There is a constant $M=M(k,\ell)$ depending on $k$,
and $\ell=\ell(\Gamma)$ such that for any given $f\in L^{k+1}_{\Z^d}(\Gamma,Y)$   and $x\in\Gamma$  there
are functions $P_x, g_x\in Y^\Gamma$ such that
for each finite subset $\Sigma\subset \Gamma$,  any $y\in\Sigma$, and any $0\leq i\leq k$,
\begin{enumerate}
\item  $f(y)=P_x(y) + g_x(y)$ \,,
\item  $P_x(y) = \sum_{m=0}^k \frac{1}{m!}\, \Theta^m  f(x)\,(y-x)^{(m)}  $ is a degree $k$ polynomial\,,
\item  $f(x) = P_x(x)$ \,,
\item  $\Theta^i (g_x) (x)=0$\,,
\item  $\spnrm{\Delta^i (g_x) }_{\Sigma}  \leq 
M\,r^{k-i+1} \,\spnrm{\Delta^{k+1} f }_\Gamma  $\, where $r=r_{x}(\Sigma)$\,,
\item  $\spnrm{\Delta^{k+1} (g_x) }_{\Gamma}  =  \spnrm{\Delta^{k+1} f }_\Gamma  $.
\end{enumerate}
\end{prop}

\bigskip

\dem
Fix any pair $(\Sigma, x)$ with  $x\in\Sigma\subset \Gamma$, and 
$\Sigma$ finite, and consider the following claim, which also
depends on the parameter $s\in\N$:

\smallskip

$\mathcal{S}(s,x)$\quad $:\Leftrightarrow$\quad
there exists $h_{s,x}\in \oplus_{i=s}^k \Hom{i}{\R^d}{Y}$
such that $g_{s,x} = f- h_{s,x}$ satisfies:
\quad\begin{enumerate}
\item[i)]\,\, $\Theta^i g_{s,x}(x)=0$, for each $i=s,\ldots, k$,
\item[ii)]\,  $\spnrm{ \Delta^i g_{s,x}}_{\Sigma}\leq 
\left( \prod_{j=i}^k
\left( \frac{j}{2}\nrm{\Gamma}+r_x(N_{j-1}\ldots N_i \Sigma) \right) \right)
\,\spnrm{\Delta^{k+1} f }_{N_k \cdots N_i \Sigma } $, 
($s\leq i\leq k$),
\item[iii)] $\spnrm{ \Delta^{k+1} g_{s,x} }_\Gamma =  \spnrm{\Delta^{k+1} f }_\Gamma$.
\end{enumerate}
We shall prove, by regressive induction, the claim
$\mathcal{S}(0,x)$.
By proposition~\ref{reducao2}, the claim
$\mathcal{S}(k,x)$ holds with
$p_{k,x}(y)=h_{k,x}(y) = \frac{1}{m!}\, \Theta^k  f(x)\,(y-x)^{(k)}$.
It lefts to show that \;
$$ \mathcal{S}(s,x)\; \Rightarrow \mathcal{S}(s-1,x),\quad
\text{for any }\; 1\leq s\leq k\;.$$
Assume that
$\mathcal{S}(s,x)$ holds and let us apply again
proposition~\ref{reducao2} to the function $g_{s,x}$.
By $\mathcal{S}(s,x)$, we have
$\spnrm{\Delta^s g_{s,x}}_\Sigma \leq
\left( \prod_{j=s}^k
\left( \frac{j}{2}\nrm{\Gamma}+r_x(N_{j-1}\ldots N_s \Sigma) \right) \right)
\,\spnrm{\Delta^{k+1} f }_{N_k \cdots N_s \Sigma }$.
Apply proposition~\ref{reducao2} to the function $g_{s,x}$ and the set
$\Sigma'=N_{s-1} \Sigma$ to  get an homogeneous polynomial
~$p_{s-1,x}\in\Hom{s-1}{\R^d}{Y}$  such that
\begin{enumerate}
\item[a)]
$\spnrm{\Delta^{s-1} g_{s-1,x}}_\Sigma \leq
\left( \prod_{j=s-1}^k
\left( \frac{j}{2}\nrm{\Gamma}+r_x(N_{j-1}\ldots N_{s-1} \Sigma) \right) \right)
\,\spnrm{\Delta^{k+1} f }_{N_k \cdots N_{s-1} \Sigma }$,\,
where
$g_{s-1,x}=g_{s,x}-p_{s-1,x}$\; and 
\item[b)] $\Theta^{s-1} g_{s-1,x}(x)=0$\;.
\end{enumerate}
We define\; $h_{s-1,x}=h_{s,x}+p_{s-1,x}$\;.
Then
\begin{eqnarray*}
g_{s-1,x} &=& g_{s,x}-p_{s-1,x} = (f-h_{s,x})- p_{s-1,x}\\
&=& f-(h_{s,x}+p_{s-1,x}) = f- h_{s-1,x}\;.
\end{eqnarray*}
By proposition~\ref{derivn},
$\Theta^i p_{s-1,x}\equiv 0$, for all $i\geq s$.
Therefore, by the induction hypothesis \, $\mathcal{S}(s,x)$\,
we have for each $i\geq s$ 
\quad\begin{enumerate}
\item[i)]\,\, $\Theta^i g_{s-1,x}(x)=\Theta^i g_{s,x}(x) +
\Theta^i p_{s-1,x}(x) = 0$, 
\item[ii)]\,  $\spnrm{\Delta^{i} g_{s-1,x}}_\Sigma =
\spnrm{\Delta^{i} g_{s,x}}_\Sigma \leq
\left( \prod_{j=i}^k
\left( \frac{j}{2}\nrm{\Gamma}+r_x(N_{j-1}\ldots N_{i} \Sigma) \right) \right)
\,\spnrm{\Delta^{k+1} f }_{N_k \cdots N_{i} \Sigma }$,
\item[iii)] $ \spnrm{ \Delta^{k+1} g_{s-1,x} }_\Gamma 
= \spnrm{\Delta^{k+1} f }_\Gamma$.
\end{enumerate}
Items a), b), i), ii) e iii) show that
$\mathcal{S}(s-1,x)$ holds.
Finally, the statement $\mathcal{S}(0,x)$ proves this proposition
with $g_x=g_{0,x}$ and $P_x=h_{0,x}=\sum_{i=0}^k p_{i,x}$.
To establish item 5. we only have to remark that when $\Sigma$ contains more than one point, 
if $r=r_x(\Sigma)$ then $\nrm{\Gamma}\leq \ell\,r$,\, and
$$ r_x(N_{j-1}\ldots N_i\Sigma) 
\leq r_x(\Sigma) +c\,\nrm{\Gamma} \leq r+c\,\ell\,r = (1+c\,\ell)\,r \;,
$$
where $c=(j-1)+(j-2)+\ldots + i$.
\cqd

\bigskip

\begin{prop}\label{Delta_m:Delta_k+1}
Given $f\in L^{k+1}_{\Z^d}(\Gamma,Y)$   and $1\leq m \leq k+1$,
$$ \spnrm{\Delta^m f }_\Gamma \leq d^{k+1-m}\,\spnrm{\Delta^{k+1} f}_\Gamma\;.$$
\end{prop}

\dem
The proof goes by induction in $k$. For $k=0$ there is nothing to prove.
Assume this inequality holds for some $k$.
Take any vector $u\in\Gamma$ and consider the zero average function
$g=\Delta_u f$. By induction hypothesis, for any vector $v \in \Gamma^m$,
\begin{align*}
\nrm{ \Delta_{(u,v)}^{m+1} f(x) } &=
\nrm{\Delta^m_v g(x)} \leq \spnrm{\Delta^{m} g}_\Gamma\,\nrm{v}^m\\
&\leq d^{k+1-m}\,\spnrm{\Delta^{k+1} g}_\Gamma\,\nrm{v}^m\\
&\leq d^{k+1-m}\,\spnrm{\Delta^{k+2} f}_\Gamma\,\nrm{u}\,\nrm{v}^m\; .
\end{align*}
Therefore\, $\spnrm{\Delta^{m+1} f}_\Gamma \leq d^{k+2-(m+1)}\,
\spnrm{\Delta^{k+2} f}_\Gamma$.
It lefts to prove the inequality for $m=1$ of the induction step $k+1$.
But we also get 
$$\spnrm{\Delta^{1} g}_\Gamma \leq d^{k}\,\spnrm{\Delta^{k+2} f}_\Gamma\,\nrm{u}\;.$$
Whence by lemma~\ref{Delta:average} 
we get \, $\nrm{g(x)} \leq d^{k+1}\,\spnrm{\Delta^{k+2} f}_\Gamma\,\nrm{u}$, \,
which implies that \,
$\spnrm{\Delta^1 f}_\Gamma \leq d^{k+1}\,\spnrm{\Delta^{k+2} f}_\Gamma$.
\cqd

We say that  
function $f\in\Fscr_{\Z^d}(\Gamma,Y)$ has
{\em zero average }\, iff\,
$\sum_{ \overline{x}\in\Gamma/\Z^d} f(\overline{x}) = 0$.

\begin{lema} \label{Delta:average}
Given $f\in L^{1}_{\Z^d}(\Gamma,Y)$   with zero average, for every $x\in\Gamma$
$$ \nrm{f(x) }  \leq d \,\spnrm{\Delta^{1} f}_\Gamma\;.$$
\end{lema}

\dem
We base the proof in the following fact. Let $D_R(y)$ denote the disk
of radius $R$ in $Y$ centered at  $y\in Y$. Then given any points
$y_1,\ldots, y_n\in Y$ we have\, 
$\cap_{i=1}^n D_R(y_i) \subseteq D_R (\,\sum_{i=1}^n y_i/n\, )$.
Let $\Gamma/\Z^d=\{\overline{x}_1,\ldots \overline{x}_n\}$ and
set $y_i=f(\overline{x}_i)$ for $i=1,\ldots, n$.
Write $C=\spnrm{\Delta^{1} f}_\Gamma$ and $R=d\, C$. We claim that
$f(\Gamma)\subseteq D_{R}(y_i)$ for each $i=1,\ldots, n$,
and the lemma follows from the previous remark.
To finish we just have to prove the claim. Given $x\in\R^d$,
we can choose a representative $x_i\in\Gamma$ of $\overline{x}_i$ 
such that $\nrm{x-x_i}\leq d$.
Then \, $\nrm{f(x)-y_i }\leq \nrm{f(x)-f(x_i)} \leq C\,\nrm{x-x_i} \leq R$, which proves that
$f(\Gamma)\subset D_R(y_i)$.
\cqd

\bigskip

\begin{prop} \label{Mtilde}
There is a constant $\widetilde{M}=\widetilde{M}(k,\ell)$ depending on $k$,
and $\ell=\ell(\Gamma)$ such that under the same assumptions of proposition~\ref{existence:px}
the polynomials there referred satisfy: for every $x,y\in\Gamma$, $u\in\Gamma^k$  and 
$0\leq m \leq k$,\,
\begin{enumerate}
\item $P_x(x)=f(x)$,
\item $ \nrm{ D^m_u (P_x)(x)}\leq\widetilde{M}\, \spnrm{\Delta^{k+1} f }_\Gamma \,\nrm{u}^m $,
\item $ \nrm{ D^m_u (P_y-P_x)(x)}\leq\widetilde{M}\,\nrm{x-y}^{k+1-m}\, \spnrm{\Delta^{k+1} f }_\Gamma \,\nrm{u}^m $.
\end{enumerate}
\end{prop}

\dem
Item 1 follows from item 3 of proposition~\ref{existence:px}.
From item 2 of the same proposition it follows that $D^m_u (P_x)(x)=\Theta^m_u f(x)$.
Whence by inequality (a) of proposition~\ref{Delta:Theta:ineq} and  proposition~\ref{Delta_m:Delta_k+1} \,
\begin{align*}
\nrm{ D^m_u (P_x)(x)} &\leq \spnrm{\Theta^m_u f}_\Gamma\, \nrm{u}^m
\leq   \spnrm{\Delta^m_u f}_\Gamma\, \nrm{u}^m \\
&\leq   d^{k+1-m} \spnrm{\Delta^{k+1}_u f}_\Gamma\, \nrm{u}^m \;.
\end{align*}
This shows that item 2. here will hold if we choose $\widetilde{M}\geq d^k$. 

It is enough to prove inequality 3. for all basic multi-vectors $u\in B^m$.
Since $P_y(z)+g_y(z)=f(z)=P_x(z)+g_x(z)$ we have 
$$\Theta^m_u( P_y-P_x)(x) =\Theta^m_u (g_x -g_y)(x) =
\underbrace{\Theta^m_u (g_x)(x) }_{=0} -\Theta^m_u (g_y)(x)= -\Theta^m_u (g_y)(x)\;.$$
Therefore, letting  $M$ be the constant of proposition~\ref{existence:px},
$\Sigma$ be a set containing $x$ and $y$ with $r=\nrm{x-y}=r_y(\Sigma)$,
by this proposition we have
\begin{align}
 \nrm{ \Theta^m_u( P_y-P_x)(x) }
& = \nrm{ \Theta^m_u (g_y)(x) } \leq \spnrm{\Delta^m (g_y)}_{N_m \Sigma}\,\nrm{u}^m\nonumber \\
& \leq M\, (r_y(N_m\Sigma) )^{k+1-m} \, \spnrm{\Delta^{k+1} f }_\Gamma \,\nrm{u}^m\nonumber \\
&\leq 
M'\,r^{k+1-m}\,\spnrm{\Delta^{k+1} f }_\Gamma \,\nrm{u}^m\label{const:M'}
\end{align}
for some universal constant $M'=M'(k,\ell)$.
Notice that 
$$r_y(N_m\Sigma)\leq r_y(\Sigma)+m\nrm{\Gamma}\leq (1+m\,\ell)\,r\;.$$
We shall now prove by regressive induction in $m$ that some
constant $\widetilde{M}_m$ exists for which inequality 3. holds
for all basic multi-vectors $u\in B^m$,
and in the end we take  $\widetilde{M}=\max_{0\leq m \leq k} \widetilde{M}_m$.
By proposition~\ref{Theta:polynomial} we have
$$ D^m_u( P_y-P_x)(x) = \Theta^m_u( P_y-P_x)(x) + R_m(x,y,u)\;,  $$
where $ R_m(x,y,u)$ denotes the remainder of  proposition~\ref{Theta:polynomial}.
For $m=k$ and $m=k-1$ the remainder vanishes, $R_m(x,y,u)=0$,
and we can take $\widetilde{M}_m = M'$, the constant in~(\ref{const:M'}).
Assume now that
$ \nrm{ D^n_u (P_y-P_x)(x)}\leq\widetilde{M}_n\,\nrm{x-y}^{k+1-n}\, \spnrm{\Delta^{k+1} f }_\Gamma \,\nrm{u}^n $
for every $n\geq m$. The remainder $R_{m-1}(x,y,u)$
is a finite sum of terms of the form
\begin{equation}
\label{Rm-1}
\frac{\tilde{r}!}{r_0!\,r_1!\cdots r_{m-1}! }\, 
D^{\tilde{r}}_{((x-y)^{(r_0)}, u_1^{(r_1)},\ldots, u_{m-1}^{(r_{m-1})} )} ( P_y-P_x)(x) \;,
\end{equation}
where $\tilde{r}\geq m+1$, $r_i$ is odd for all $1\leq i\leq m-1$, $r_i>1$ for some $1\leq i\leq m-1$
 and $r_0+r_1+\ldots + r_{m-1}=\tilde{r}$.
Assume $x\neq y$. Otherwise inequality 3. is clear.
Now, since $u$ is a basic multi-vector we must have \,
$\nrm{u_i}\leq \ell \,\nrm{x-y}$,\, which implies
$\nrm{u_i}^{r_i}\leq \ell^{r_i-1} \,\nrm{x-y}^{r_i-1}\, \nrm{u_i}$.
By induction hypothesis term~(\ref{Rm-1}) is bounded by
\begin{align}
& \frac{\tilde{r}!}{r_0!\,r_1!\cdots r_{m-1}! }\, 
\widetilde{M}_{\tilde{r}}\,\nrm{x-y}^{k+1-(r_1+\cdots+r_{m-1})}\, \spnrm{\Delta^{k+1} f }_\Gamma\,
\nrm{u_1}^{r_1}\cdots \nrm{u_{m-1}}^{r_{m-1}} \nonumber \\
&\leq   \frac{\tilde{r}!}{r_0!\,r_1!\cdots r_{m-1}! }\, 
\widetilde{M}_{\tilde{r}}\,\ell^{r_1+\cdots+r_{m-1}-(m-1)}\, \nrm{x-y}^{k+1-(m-1)}\, \spnrm{\Delta^{k+1} f }_\Gamma\,
\nrm{u}^{m-1} \label{term:bound} \;.
\end{align}
Define
$$\widetilde{M}_{m-1} := M' + \sum_{r_0+r_1+\cdots + r_{m-1}=\tilde{r}}
 \frac{\tilde{r}!}{r_0!\,r_1!\cdots r_{m-1}! }\, 
\widetilde{M}_{\tilde{r}}\,\ell^{r_1+\cdots+r_{m-1}-(m-1)}\;, $$
the sum being taken over the same set of multi-indices 
we took in~(\ref{Rm-1}).
Combining~(\ref{const:M'}) (for  $m-1$) with~(\ref{term:bound}),
we see that
$$ \nrm{ D^{m-1}_u (P_y-P_x)(x)}\leq\widetilde{M}_{m-1}\,\nrm{x-y}^{k+1-({m-1})}\, \spnrm{\Delta^{k+1} f }_\Gamma \,\nrm{u}^{m-1} \;,$$
and this completes the induction step. 
\cqd

\bigskip

%

\demof{Proof of Theorem A}
Consider the constant $A=A(d,k)$ provided by Fefferman's theorem
(stated in the introduction) and set\, 
$\widetilde{A}(d,k,\ell):=\widetilde{M}(k-1,\ell)\,A(d,k)$,
where $\widetilde{M}=\widetilde{M}(k-1,\ell)$ is the constant given by
proposition~\ref{Mtilde}. Any function $f:\Gamma/\Z^d\to Y$ with
bounded differences of order $k$ lifts to a function  $f\in L^{k}_{\Z^d}(\Gamma,Y)$.
Take any exhausting family of finite subsets
$\Sigma_n\subset \Gamma$,\,  $\Gamma =\bigcup_{n=0}^\infty \Sigma_n$,
\, $\Sigma_1\subset \Sigma_2\subset \ldots\Sigma_n \subset \ldots \,$.
By proposition~\ref{Mtilde} the function
$f\vert_{\Sigma_n}:\Sigma_n\to Y$ satisfies the $C^k$-Whitney extension condition
with constant $\widetilde{M}\,\spnrm{\Delta^k f}_\Gamma$,
where $\widetilde{M}=\widetilde{M}(k-1,\ell)$.
By Fefferman's theorem there is a function $F_n:\R^d\to Y$
such that $F_n(x)=f(x)$ for every
$x\in   \Sigma_n$, and 
$\spnrm{D^k F_n}_{\R^d} \leq A \widetilde{M} \,\spnrm{\Delta^k f}_{\Gamma} = 
\widetilde{A} \,\spnrm{\Delta^k f}_{\Gamma} $.
Consider now the convex set
$$\Fscr=\{\, F\in\CD{k-1,1}{\R^d}{Y}\,:
\, {\rm Lip}(D^{k-1} F) \leq \widetilde{A} \,\spnrm{\Delta^k f}_{\Gamma}\;\} \;. $$
The set $\Fscr$ is compact for the topology of uniform convergence 
(of all derivatives up to order $k-1$) over compact sets in $\R^d$.
Since $F_n\in\Fscr$ for each $n\geq 0$, there is a subsequence
$F_{n_j}$ which converges to a function $F\in\Fscr$.
Because $F_n(x) = f(x)$ for all $n\geq p$ and $x\in\Sigma_p$,
it follows that $F(x) = f(x)$ for all  $x\in\Gamma$.
This proves that
$$\Fscr_f=\{\, F\in\Fscr\,:\, F(x) = f(x),\;\forall x\in\Gamma\;\}$$
is a non-empty compact convex subset of $\Fscr$.
For each vector $v\in\R^d$ consider the translation operator
$\tau_v:\CD{k-1,1}{\R^d}{Y}\to \CD{k-1,1}{\R^d}{Y}$,\,
$\tau_v(F)(x)=F(x+v)$ ($x\in \R^d$). The set $\Fscr_f$ is invariant
under all translations $\tau_v$ with $v\in\Z^d$, because  $f$ is $\Z^d$-periodic.
Given vectors $v_1,\ldots, v_i\in\Z^d$ we denote by $\Fscr_f[v_1,\ldots, v_i]$
the set of all functions $F\in\Fscr_f$ such that
$F=\tau_{v_1} F = \ldots = \tau_{v_i} F $.
Assume now that $\{v_1,\ldots, v_d\}$ is a basis for $\Z^d$.
Then $\Fscr_f[v_1,\ldots, v_d]$ is the set of $\Z^d$-periodic extensions
of $f$ in $\Fscr_f$. Our goal is to prove that $\Fscr_f[v_1,\ldots, v_d]\neq \emptyset$,
which is done by induction in the number $i$ of vectors in $\Fscr_f[v_1,\ldots, v_i]$.
For $i=0$, $\Fscr_f[\,]=\Fscr_f$ and we already know that $\Fscr_f\neq\emptyset$.
Assume $\Fscr_f[v_1,\ldots, v_{i-1}]\neq \emptyset$ and take $F\in \Fscr_f[v_1,\ldots, v_{i-1}]$.
We notice that the set  $\Fscr_f[v_1,\ldots, v_{i-1}]$ is always compact, convex and invariant
under all translations $\tau_v$ with $v\in\Z^d$.
Then  the functions  $S_n= \frac{1}{n}\sum_{j=0}^{n-1} \tau_{j\,v_{i}}F$
belong to $\Fscr_f[v_1,\ldots, v_{i-1}]$, and because this is a compact set,
there is a  subsequence $S_{n_j}$ which converges to some function 
$\widetilde{F}\in \Fscr_f[v_1,\ldots, v_{i-1}]$.
Finally, since $S_n-\tau_{v_i} S_n =\frac{1}{n}\left(F-\tau_{n v_i} F \right)$
converges to zero, it follows that $\widetilde{F}\in \Fscr_f[v_1,\ldots, v_{i}]$.
\cqd

\bigskip

\thispagestyle{empty}

\end{document}